# Spectral properties of Kac-Murdock-Szegö matrices with a complex parameter

George Fikioris

School of Electrical and Computer Engineering

National Technical University of Athens

GR 157-73 Zografou, Athens, Greece

email: gfiki@ece.ntua.gr

*Abstract*

When $0 < \rho < 1$, the Kac-Murdock-Szegö matrix $K_n(\rho) = \left[\rho^{|j-k|}\right]_{j,k=1}^{n}$ is a Toeplitz correlation matrix with many applications and very well known spectral properties. We study the eigenvalues and eigenvectors of $K_n(\rho)$ for the general case where $\rho$ is complex, pointing out similarities and differences to the case $0 < \rho < 1$. We then specialize our results to real $\rho$ with $\rho > 1$, emphasizing the continuity of the eigenvalues as functions of $\rho$. For $\rho > 1$, we develop simple approximate formulas for the eigenvalues and pinpoint all eigenvalues' locations. Our study starts from a certain polynomial whose zeros are connected to the eigenvalues by elementary formulas. We discuss relations of our results to earlier results of W. F. Trench.

## 1. Introduction

For $n \geq 2$ and $\rho \in \mathbf{C}$, let $K_n(\rho)$ be the $n \times n$ complex symmetric Toeplitz matrix

$$K_n(\rho) = \left[\rho^{|j-k|}\right]_{j,k=1}^{n} = \begin{bmatrix} 1 & \rho & \rho^2 & \ldots & \rho^{n-1} \\ \rho & 1 & \rho & \ldots & \rho^{n-2} \\ \rho^2 & \rho & 1 & \ldots & \rho^{n-3} \\ \vdots & \vdots & \vdots & \ddots & \vdots \\ \rho^{n-1} & \rho^{n-2} & \rho^{n-3} & \ldots & 1 \end{bmatrix}. \qquad (1.1)$$

In the special case $0 < \rho < 1$, the real symmetric and positive definite matrix $K_n(\rho)$ is often called the Kac-Murdock-Szegö matrix, the term originating from the 1953 paper [1]. In [1] (see also the book [2] by Grenander and Szegö), a detailed study of $K_n(\rho)$ was an important first step for the development of general theorems on the asymptotic (large-$n$) behavior of eigenvalues and determinants of Toeplitz matrices. Furthermore, [1] and [2] illustrated their theorems by specialization to $K_n(\rho)$.

When $0 < \rho < 1$, $K_n(\rho)$ is a correlation matrix whose spectral properties are very well known [1]-[4] and whose applications are pervasive. It is important within the context of digital signal processing because it characterizes first-order stationary Markov random signals, which can be decorrelated via the Karhunen-Loéve Transform. Thus $K_n(\rho)$ plays a central role in the theory of this transform, as well as related transforms such as the Discrete Cosine Transform [5]-[8]. As the autocorrelation matrix of stationary AR(1) processes, $K_n(\rho)$ is important for models of multidimensional scaling which give rise to horseshoes (quadratic curves) [4], [9]. It is used in [10], which studies a certain class of matrices related to principal component analysis and the yield curve of multifactor models of interest rates. It is the "exponential correlation matrix" of the work [11] pertaining to MIMO (multiple-input multiple output) channels in wireless communications. The matrix $K_n(\rho)$ is a test matrix in MATLAB [12], one of the "Famous Matrices" of GNU Octave [13], and is often used as a means of illustration of general theorems—e.g. for $0 < \rho < 1$ $K_n(\rho)$ is one of the example matrices of [14] (see also [15]), where Bogoya, Böttcher, Grudsky, and Maximenko develop asymptotics for the eigenvectors of Hermitian Toeplitz matrices. When $\rho \in \mathbf{C}$ with $|\rho| < 1$, $K_n(\rho)$ characterizes one-dimensional periodic lattices with geometric interaction, a problem relevant to quantum emitters acting as quantum information processing units and, also, to discrete fiber gratings used for frequency selection [16].

The purpose of the present paper is to go beyond the well-known case $0 < \rho < 1$. Sections 2-5 are a standalone discussion of the eigensystem of $K_n(\rho)$ for the general case $\rho \in \mathbf{C}$. We show that many facts, known for the case $0 < \rho < 1$, carry

over—often with modifications—to the general case. Sections 6 and 7 specialize to $\rho \in \mathbf{R}$ with $\rho > 1$; in this case, much more can be said.

As opposed to the case $0 < \rho < 1$, $K_n(\rho)$ can admit double eigenvalues for certain values of $\rho \in \mathbf{C}$, as shown in Section 4. The cases $0 < \rho < 1$ and $\rho > 1$ are also dissimilar. An essential difference follows from the more general study of Trench [17] (see also Trench's [18], as well as the notes by Bogoya, Böttcher, and Grudsky in [19, p. 275]): When $\rho > 1$, and for sufficienly large $n$, it follows from [17] that $K_n(\rho)$ possesses certain "wild" or "outlying" or "un-Szegö-like" eigenvalues. In Section 6, we shed further light on these issues by developing precise conditions on $\rho$ and $n$ for the existence of what we call "extraordinary eigenvalues" of $K_n(\rho)$ and by investigating the behavior of these eigenvalues in detail.

When $|\rho| > 1$, $K_n(\rho)$ no longer belongs to certain matrix classes that facilitate studies of its spectral behavior. For example $K_n(\rho)$ is not positive definite. Perhaps more importantly, the associated Laurent matrix does not have a well-defined and bounded symbol. Thus an exploration of the spectrum for the general case $\rho \in \mathbf{C}$ requires other methods of study. In this paper, our initial tool is a polynomial denoted by $p_{2n}(\rho, z)$. This polynomial (which is self-inversive in the case $\rho \in \mathbf{R}$) is defined by elementary closed-form expressions and its zeros are connected to the eigenvalues of $K_n(\rho)$ by simple formulas.

We begin in Section 2 by giving some known formulas related to $K_n(\rho)$, e.g. for its characteristic polynomial and inverse. We then provide necessary and sufficient conditions for $K_n(\rho)$ to belong to certain matrix classes that have proved useful for the case $0 < \rho < 1$. Section 3 discusses $p_{2n}(\rho, z)$ and explains the significance of zeros that lie on the unit circle. In Section 4 we show that $p_{2n}(\rho, z)$ naturally gives rise to trigonometric functions that have, in the literature, been utilized for the determination of the eigensystem when $\rho \in \mathbf{R}$, thus demonstrating that the same trigonometric functions are useful, more generally, for $\rho \in \mathbf{C}$. We subsequently use these functions to investigate the issue of double eigenvalues. For $\rho \in \mathbf{C}$ with $|\rho| > 1$,

Section 5 demonstrates that $K_n(\rho)$ has two "large" eigenvalues and discusses their relation with results of Trench. For the case $\rho \in \mathbf{R}$ with $\rho > 1$, Sections 6 and 7 analyze in detail the eigensystem of $K_n(\rho)$; included here are bounds and simple approximations for the eigenvalues, as well as a number of different formulas useful for numerical computations. In particular, these two sections emphasize the continuity of the $n$ eigenvalues as functions of the real parameter $\rho$; to this end, we are helped by repeating some results from the literature that pertain to the case $0 < \rho < 1$.

Throughout, we assume $n = 2, 3, \ldots$. When counting matrix eigenvalues and polynomial zeros, we always take multiplicities into account. A notation such as $p_{2n}(\rho, z)$ denotes a monic polynomial in $z$ of degree $2n$ with $\rho$-dependent coefficients. When a polynomial is given by an expression like (3.1) below, it will go without saying that all denominator zeros ($z = 1$ and $z = -1$ in the case of (3.1)) are removable singularities. Similarly, (3.8) and (3.9) together constitute a formula for $p_{2n}(\rho, z)$ in which $z = \rho$ and $z = 1/\rho$ are, necessarily, removable singularities.

A Toeplitz matrix $\left[a_{j-k}\right]_{j,k=1}^{n}$ is often considered to be a truncation (or finite section) of the Laurent or infinite-Toeplitz matrices respectively given by

$$\begin{bmatrix} \ldots & \ldots & \ldots & \ldots & \ldots \\ \ldots & a_0 & a_{-1} & a_{-2} & \ldots \\ \ldots & a_1 & a_0 & a_{-1} & \ldots \\ \ldots & a_2 & a_1 & a_0 & \ldots \\ \ldots & \ldots & \ldots & \ldots & \ldots \end{bmatrix}, \quad \begin{bmatrix} a_0 & a_{-1} & a_{-2} & \ldots & \ldots \\ a_1 & a_0 & a_{-1} & \ldots & \ldots \\ a_2 & a_1 & a_0 & \ldots & \ldots \\ \ldots & \ldots & \ldots & \ddots & \ldots \\ \ldots & \ldots & \ldots & \ldots & \ldots \end{bmatrix}. \qquad (1.2)$$

When the Fourier series $\sum_{k=-\infty}^{\infty} a_k e^{ik\theta}$ ($\theta \in \mathbf{R}$) is a well-defined function $\sigma(\theta)$, then $\sigma(\theta)$ is called the *symbol* (see, e.g. [20], [21]). If $\sigma(\theta)$ is bounded and continuous, then the spectrum of the Laurent matrix in (1.2) is the *range* of $\sigma(\theta)$, i.e., the continuous set of complex values assumed by $\sigma(\theta)$ for $\theta \in (-\pi, \pi]$. We will denote this set by $\text{range}\{\sigma(\theta)\}$.

We conclude this Introduction with some facts about self-inversive polynomials [22], [23]. A classical theorem due to Cohn states that all zeros of the $m$'th degree polynomial

$$f(z) = \sum_{k=0}^{m} a_k z^k, \quad a_k \in \mathbf{C}, \quad a_m \neq 0 \tag{1.3}$$

lie on the unit circle iff (i) $f(z)$ is self-inversive; and (ii) all zeros of the derivative $f'(z)$ lie in or on the unit circle. The polynomial $f(z)$ is called *self-inversive* if

$$a_k = \varepsilon \bar{a}_{m-k}, \quad k = 0, 1, \ldots, m, \quad \varepsilon \in \mathbf{C}, \quad |\varepsilon| = 1, \tag{1.4}$$

where the overbar denotes the complex conjugate. When all coefficients $a_0, a_1, \ldots, a_m$ are real, $f(z)$ can be self-inversive only when $\varepsilon = \pm 1$. In the first case ($\varepsilon = 1$), the self-inversive polynomial is called *reciprocal*; when $\varepsilon = -1$, the self-inversive polynomial is called *anti-reciprocal*.

## 2. Formulas and matrix classes related to $K_n(\rho)$

We begin by giving some known properties of $K_n(\rho)$ that will prove useful to us.

**Theorem 2.1** Let $\rho \in \mathbf{C}$.

(i) For $\rho \neq \pm 1$, the *inverse* $[K_n(\rho)]^{-1}$ is the non-Toeplitz tridiagonal matrix given by

$$[K_n(\rho)]^{-1} = \frac{1}{1-\rho^2} \begin{bmatrix} 1 & -\rho & 0 & \cdots & 0 & 0 & 0 \\ -\rho & 1+\rho^2 & -\rho & \cdots & 0 & 0 & 0 \\ 0 & -\rho & 1+\rho^2 & \cdots & 0 & 0 & 0 \\ \vdots & \vdots & \vdots & \ddots & \vdots & \vdots & \vdots \\ 0 & 0 & 0 & \cdots & 1+\rho^2 & -\rho & 0 \\ 0 & 0 & 0 & \cdots & -\rho & 1+\rho^2 & -\rho \\ 0 & 0 & 0 & \cdots & 0 & -\rho & 1 \end{bmatrix}. \tag{2.1}$$

(ii) The *characteristic polynomial* $\psi_n(\rho, \lambda) = \det[\lambda I_n - K_n(\rho)]$ satisfies the recurrence relation

$$\psi_n(\rho, \lambda) = [\rho^2 - 1 + \lambda(1+\rho^2)]\psi_{n-1}(\rho, \lambda) - (\lambda \rho)^2 \psi_{n-2}(\rho, \lambda), \tag{2.2}$$

whose appropriate initial conditions are

$$\psi_0(\rho,\lambda)=1, \quad \psi_1(\rho,\lambda)=\lambda-1. \tag{2.3}$$

The characteristic polynomial is explicitly given by

$$\psi_n(\rho,\lambda)=\frac{(\lambda\rho)^n}{1-\rho^2}\left[U_n(\tau)-2\rho U_{n-1}(\tau)+\rho^2 U_{n-2}(\tau)\right], \quad n=2,3,\ldots, \tag{2.4}$$

where

$$\tau=\tau(\rho,\lambda)=\frac{\rho^2(\lambda+1)+\lambda-1}{2\lambda\rho}, \tag{2.5}$$

and where $U_k(t)$ is the Chebyshev polynomial of the second kind, defined by

$$U_k(t)=\frac{\sin[(k+1)\theta]}{\sin\theta}, \quad t=\cos\theta\in\mathbf{C}, \quad k=0,1,2,\ldots. \tag{2.6}$$

For $\rho=0$ and $\rho=\pm 1$, in particular, we have

$$\psi_n(0,\lambda)=(\lambda-1)^n, \quad \psi_n(\pm 1,\lambda)=\lambda^{n-1}(\lambda-n), \quad n=1,2,\ldots. \tag{2.7}$$

(iv) The *determinant* of $K_n(\rho)$ equals $(1-\rho^2)^{n-1}$ and $K_n(\rho)$ is singular iff $\rho=\pm 1$.

(v) The $\rho$-dependent *symbol* $\sigma(\rho,\theta)$ of the corresponding Laurent matrix is well defined and bounded if and only if $\rho\in\mathbf{C}$ with $|\rho|<1$, in which case it equals the Poisson kernel for the unit disc, viz.,

$$\sigma(\rho,\theta)=\frac{1-\rho^2}{1-2\rho\cos\theta+\rho^2}, \quad \theta\in(-\pi,\pi]. \tag{2.8}$$

*Proof:* Eqn. (2.1), which can be verified directly, is in [24, 7.2.P13, p. 443] and [3]. Eqns. (2.2)-(2.6) are shown in [1] and [2] (these two works assume $0<\rho<1$, but their proof holds for all $\rho\in\mathbf{C}$). When $\rho=\pm 1$ or $\rho=0$, (2.4) with (2.5) become indeterminate, with their limiting values given by (2.7); eqns. (2.7) can be verified using (2.2) and (2.3). The determinant is in [24, 7.2.P12, p. 443]. Eqn. (2.8), which can be found in [1]-[3], follows by directly summing the defining Fourier series. □

**Remark 2.2** The right-hand side of (2.8) has a straightforward analytic continuation to $|\rho| \geq 1$, and we will continue to denote this analytic continuation by $\sigma(\rho, \theta)$. Similarly, range$\{\sigma(\rho, \theta)\}$ will denote the range of *both* the quantity in (2.8) (i.e., the symbol of the Laurent matrix) *and* its analytic continuation. While the analytic continuation should not be confused with the symbol of some Laurent matrix, we will see that range$\{\sigma(\rho, \theta)\}$ is relevant to the finite matrix $K_n(\rho)$ both when $|\rho| < 1$ and when $|\rho| \geq 1$.

We now give conditions for $K_n(\rho)$ to be in a number of matrix classes. As already mentioned in our Introduction, satisfaction of one (or more) of these conditions facilitates a spectral analysis of $K_n(\rho)$.

**Theorem 2.3** For $\rho \in \mathbf{C}$, Table 1 gives necessary and sufficient conditions in order for $K_n(\rho)$ to belong to certain matrix classes:

|       | **Property** | **Necessary and sufficient conditions** |
|-------|--------------|------------------------------------------|
| (i)   | $K_n(\rho)$ is positive | $\rho > 0$ |
| (ii)  | $K_n(\rho)$ is real symmetric | $\rho \in \mathbf{R}$ |
| (iii) | $K_n(\rho)$ is Hermitian | $\rho \in \mathbf{R}$ |
| (iv)  | The corresponding to $K_n(\rho)$ Laurent matrix has a well-defined and bounded symbol $\sigma(\rho, \theta)$ | $\rho \in \mathbf{C}$ with $|\rho| < 1$ |
| (v)   | $K_n(\rho)$ is positive definite (positive semidefinite) | $-1 < \rho < 1$ ($-1 \leq \rho \leq 1$) |
| (vi)  | $K_n(\rho)$ is normal | $\rho \in \mathbf{R}$ or $n = 2$ |
| (vii) | $K_n(\rho)$ is a Green's matrix | $\rho \in \mathbf{R} \setminus \{0\}$ |
| (viii)| $K_n(\rho)$ is totally positive (i.e., every minor is nonnegative, see [25, p. 11]) | $0 \leq \rho \leq 1$ |
| (ix)  | $K_n(\rho)$ is oscillatory | $0 < \rho < 1$ |

Table 1: Properties of $K_n(\rho)$ and corresponding necessary and sufficient conditions on $\rho$ and $n$

*Proof of Theorem 2.3:* (i)-(iii) are obvious while (iv) repeats Theorem 2.1(v).

(v): When $-1 < \rho < 1$, $K_n(\rho)$ is positive definite (p-d) by [24, 7.2.P12, p. 443]. Conversely, if $K_n(\rho)$ is p-d then it is Hermitian [24, p. 231] and thus $\rho \in \mathbf{R}$ by (ii). When $\rho = \pm 1$ zero is an eigenvalue by (2.7), and when $\rho > 1$ or $\rho < -1$ we will see below (Theorem 6.6 and Lemma 2.5) that $K_n(\rho)$ has at least one negative eigenvalue. Since a p-d matrix has positive eigenvalues, we are left with $-1 < \rho < 1$.

(vi): The complex matrix $K_n(\rho)$ is normal iff

$$K_n(\rho) K_n(\bar{\rho}) = K_n(\bar{\rho}) K_n(\rho) \tag{2.9}$$

Eqn. (2.9) is satisfied when ($\rho \in \mathbf{R}$ or $n = 2$), so this condition is sufficient for normality. A quick calculation gives the 1,2 element of the left-hand side of (2.9) as

$$\left[ K_n(\rho) K_n(\bar{\rho}) \right]_{12} = 2\Re\rho + \rho \sum_{k=3}^{n} |\rho|^{2(k-2)}$$

in which the sum is empty when $n = 2$. Consequently, (2.9) implies

$$2\Re\rho + \rho \sum_{k=3}^{n} |\rho|^{2(k-2)} = 2\Re\rho + \bar{\rho} \sum_{k=3}^{n} |\rho|^{2(k-2)}$$

which in turn implies ($\rho \in \mathbf{R}$ or $n = 2$). Thus our sufficient condition is necessary.

(vii) (cf. [4], [10]): By definition [25, p. 110], a real $n \times n$ matrix $A$ is a Green's matrix if its elements $A_{jk}$ can be expressed in the form

$$A_{jk} = \begin{cases} \alpha_j \beta_k, & 1 \leq j \leq k \leq n \\ \alpha_k \beta_j, & 1 \leq k \leq j \leq n \end{cases} \tag{2.10}$$

where $\alpha_j, \beta_j \in \mathbf{R}$. For $\rho = 0$ (2.10) cannot be satisfied. For $\rho \in \mathbf{R} \setminus \{0\}$, taking $\alpha_j = \beta_j^{-1} = \rho^{-j}$ shows that $K_n(\rho)$ is a Green's matrix.

(viii) (cf. [10]): Since any totally positive (TP) matrix is real, take $\rho \in \mathbf{R}$. If $\rho \in \mathbf{R} \setminus \{0\}$, then $K_n(\rho)$ is a Green's matrix by (vii). By [25, p. 111], a Green's matrix is TP iff all $\alpha_j$ and $\beta_j$ in (2.10) have the same strict sign with

$$\frac{\alpha_1}{\beta_1} \leq \frac{\alpha_2}{\beta_2} \leq \cdots \leq \frac{\alpha_n}{\beta_n}.$$

In the case $A = K_n(\rho)$, the aforesaid conditions are equivalent to $0 < \rho \leq 1$. (For $\rho = 0$, $K_n(0) = I_n$ is TP by definition.)

(ix) By [26] (but note that [26] uses the term "totally nonnegative" for our "totally positive"), a totally positive $n \times n$ matrix $A$ with elements $A_{jk}$ is oscillatory iff it is nonsingular with $A_{j,j+1} > 0$ and $A_{j+1,j} > 0$ for $j = 1, 2, \ldots, n-1$. Accordingly, the desired equivalence follows from (viii), Theorem 2.1(iv), and $K_n(0) = I_n$. □

**Remark 2.4** For $\rho \in \mathbf{C}$ with $|\rho| < 1$, Theorem 2.3(iv) states that $K_n(\rho)$ has a well-defined and bounded symbol. By a theorem in [20, p. 22-4], matrices with symbols of this type are normal for all $n = 1, 2, \ldots$ iff range$\{\sigma(\rho, \theta)\}$ is a straight line. For $\rho \in \mathbf{C}$ with $|\rho| < 1$, it is shown in [16] that range$\{\sigma(\rho, \theta)\}$ is a circular arc, and the formulas of [16] readily imply that this arc reduces to a straight line iff $-1 < \rho < 1$. This observation agrees with Theorem 2.3(vi).

When discussing the eigensystem of $K_n(\rho)$, the following lemma allows us to assume that $\rho$ lies in the first quadrant.

**Lemma 2.5** Let $(\lambda, \mathbf{y})$ be an eigenpair of $K_n(\rho)$ and let $J_n$ be the signature matrix

$$J_n = \begin{bmatrix} 1 & 0 & \cdots & 0 & 0 \\ 0 & -1 & \cdots & 0 & 0 \\ \vdots & \vdots & \ddots & 0 & 0 \\ 0 & 0 & \cdots & (-1)^n & 0 \\ 0 & 0 & & 0 & (-1)^{n+1} \end{bmatrix} \quad (2.11)$$

Then $K_n(\bar{\rho})$ and $K_n(-\rho)$ possess the eigenpairs $(\bar{\lambda}, \bar{\mathbf{y}})$ and $(\lambda, J_n \mathbf{y})$, respectively.

*Proof:* The assertion about $K_n(\bar{\rho})$ follows immediately from $\overline{K_n(\rho)} = K_n(\bar{\rho})$. Use $J_n J_n = I_n$ to rewrite $K_n(\rho)\mathbf{y} = \lambda \mathbf{y}$ as $K_n(\rho) J_n J_n \mathbf{y} = \lambda \mathbf{y}$. Left multiplication by $J_n$ yields $[J_n K_n(\rho) J_n](J_n \mathbf{y}) = \lambda (J_n \mathbf{y})$ or, by (1.1) and (2.11), $K_n(-\rho)(J_n \mathbf{y}) = \lambda (J_n \mathbf{y})$. □

## 3. Polynomials associated with eigenvalues; zeros on unit circle

For $n = 2, 3, \ldots$, $\rho \in \mathbf{C}$, and $z \neq \pm 1$ (a restriction to be removed momentarily), define

$$p_{2n}(\rho, z) = \frac{z^{2n}(z-\rho)^2 - (\rho z - 1)^2}{z^2 - 1}. \tag{3.1}$$

Evidently,

$$p_{2n}(\rho, z) = \frac{s_{n+1}(\rho, z) c_{n+1}(\rho, z)}{z^2 - 1}, \tag{3.2}$$

where

$$s_{n+1}(\rho, z) = z^{n+1} - \rho z^n + \rho z - 1 \tag{3.3}$$

$$c_{n+1}(\rho, z) = z^{n+1} - \rho z^n - \rho z + 1, \tag{3.4}$$

are monic polynomials. The expression

$$p_{2n}(\rho, z) = z^{2n} + (1+\rho^2) \sum_{k=1}^{n-1} z^{2k} - 2\rho \sum_{k=0}^{n-1} z^{2k+1} + 1, \tag{3.5}$$

which can be verified from (3.1), removes the restriction $z \neq \pm 1$ and makes it plain that, as a function of $z$, $p_{2n}(\rho, z)$ is a monic polynomial of degree $2n$ that satisfies

$$z^{2n} p_{2n}\left(\rho, \frac{1}{z}\right) = p_{2n}(\rho, z). \tag{3.6}$$

The polynomials of (3.3) and (3.4) satisfy the similar relations

$$z^{n+1} s_{n+1}\left(\rho, \frac{1}{z}\right) = -s_{n+1}(\rho, z), \quad z^{n+1} c_{n+1}\left(\rho, \frac{1}{z}\right) = c_{n+1}(\rho, z). \tag{3.7}$$

Eqn. (3.4) gives the coefficients of $c_{n+1}(\rho, z)$ as $a_0 = a_{n+1} = 1$, $a_1 = a_n = -\rho$, and $a_2 = a_3 = \ldots = a_{n-1} = 0$, where we use the notation of (1.3). Thus (1.4) is satisfied iff $\rho \in \mathbf{R}$ and $\varepsilon = 1$, meaning (see Introduction) that $c_{n+1}(\rho, z)$ is reciprocal. Reasoning similarly with (3.3) and (3.5) leads to **Lemma 3.1** The polynomials $p_{2n}(\rho,z)$, $s_{n+1}(\rho,z)$, and $c_{n+1}(\rho,z)$ are self-inversive iff $\rho \in \mathbf{R}$, in which case $p_{2n}(\rho,z)$ and $c_{n+1}(\rho,z)$ are reciprocal, while $s_{n+1}(\rho,z)$ is anti-reciprocal.

Whether self-inversive or not, $p_{2n}(\rho,z)$ is connected to the characteristic polynomial $\psi_n(\rho,\lambda)$:

**Lemma 3.2** For $\rho \in \mathbf{C}\setminus\{-1,0,1\}$, $p_{2n}(\rho,z)$ and $\psi_n(\rho,\lambda)$ are related via

$$p_{2n}(\rho,z) = \frac{(z-\rho)^n(1-z\rho)^n}{\rho^n(1-\rho^2)^{n-1}}\psi_n(\rho,\lambda), \tag{3.8}$$

in which $\lambda = \sigma(\rho,-i\ln z)$ (see Remark 2.2) or

$$\lambda = \lambda(\rho,z) = \frac{z(1-\rho^2)}{(z-\rho)(1-z\rho)}. \tag{3.9}$$

*Proof:* When $\lambda = \lambda(\rho,z)$ is given by (3.9), the $\tau$ of (2.5) becomes $\tau = (z^2+1)/(2z)$. Accordingly (2.4) yields

$$\psi_n(\rho,\lambda(\rho,z)) = \frac{z^n\rho^n(1-\rho^2)^{n-1}}{(z-\rho)^n(1-z\rho)^n}\left[U_n\left(\frac{z^2+1}{2z}\right) - 2\rho U_{n-1}\left(\frac{z^2+1}{2z}\right) + \rho^2 U_{n-2}\left(\frac{z^2+1}{2z}\right)\right]. \tag{3.10}$$

Setting $z = e^{i\theta}$ in (2.6) gives the identity

$$U_k\left(\frac{z^2+1}{2z}\right) = \frac{z^{k+2}-z^{-k}}{z^2-1}, \quad k=0,1,\ldots, \quad z \in \mathbf{C}\setminus\{0\}. \tag{3.11}$$

By (3.11) and (3.1), the quantity in square brackets in (3.10) equals $p_{2n}(\rho,z)/z^n$ and we are led to (3.8). Comparison of (2.8) and (3.9) gives $\lambda = \sigma(\rho,-i\ln z)$. □

Counting multiplicities, Lemma 3.2 provides a two-to-one correspondence between the zeros of $p_{2n}(\rho,z)$ and the eigenvalues of $K_n(\rho)$:

**Theorem 3.3** Let $\rho \in \mathbf{C}\setminus\{-1,0,1\}$. If $z_k$ is a zero of $p_{2n}(\rho,z)$, then $z_k \neq 0$ and $z_k^{-1}$ is also a zero of $p_{2n}(\rho,z)$. Furthermore

$$\lambda_k = \frac{z_k(1-\rho^2)}{(z_k-\rho)(1-\rho z_k)} = \frac{z_k^{-1}(1-\rho^2)}{(z_k^{-1}-\rho)(1-\rho z_k^{-1})} \qquad (3.12)$$

is an eigenvalue of $K_n(\rho)$. Conversely, if $\lambda_k$ is an eigenvalue of $K_n(\rho)$, then $\lambda_k \neq 0$ and the two quantities

$$z_k = \tau_k + (\tau_k^2 - 1)^{1/2}, \quad z_k^{-1} = \tau_k - (\tau_k^2 - 1)^{1/2} \qquad (3.13)$$

are zeros of $p_{2n}(\rho,z)$, where $\tau_k = \tau(\rho,\lambda_k)$ is found from (2.5).

*Proof:* Let $p_{2n}(\rho, z_k) = 0$. We have $z_k \neq 0$ by (3.5), so (3.6) yields $p_{2n}(\rho, z_k^{-1}) = 0$. Eqns. (3.1), (3.6), and the condition $\rho \neq \pm 1$ give

$$p_{2n}(\rho,\rho) = \rho^{2n} p_{2n}(\rho, 1/\rho) = 1 - \rho^2 \neq 0. \qquad (3.14)$$

Eqn. (3.14) implies $z_k \neq \rho$ and $z_k \neq 1/\rho$. Thus the $\lambda_k$ of (3.12) is well defined and comparison with (3.8) and (3.9) shows that $\psi_n(\rho, \lambda_k) = 0$. Since $\psi_n(\rho, \lambda)$ is the characteristic polynomial, we have shown (3.12).

Conversely, let $\psi_n(\rho, \lambda_k) = 0$. The eigenvalue $\lambda_k$ is nonzero by $\rho \neq \pm 1$ and Theorem 2.2(iv). Eqn. (3.9) and $\rho \neq 0$ show that $\lambda = \lambda_k$ corresponds to two values of $z$. Solving (3.9) for $z$ gives the two quantities in (3.13), with $\tau_k = \tau(\rho, \lambda_k)$ found from (2.5) (since $\rho \neq 0$ and $\lambda_k \neq 0$, $\tau_k$ is well defined). Eqn. (3.8) and $\psi_n(\rho, \lambda_k) = 0$ then entail $p_{2n}(\rho, z_k) = 0$. Finally, $p_{2n}(\rho, z_k^{-1}) = 0$ follows from (3.6). □

Determining the spectrum of $K_n(\rho)$ thus reduces to finding the zeros of the polynomial $p_{2n}(\rho,z)$, or the related through (3.2) polynomials $s_{n+1}(\rho,z)$ and $c_{n+1}(\rho,z)$. It is important that all three polynomials are given by simple formulas (e.g. $s_{n+1}(\rho,z)$ and $c_{n+1}(\rho,z)$ are lacunary in the sense discussed in [27]), and exploiting this simplicity is a major theme of this paper. To give an elementary example, for the special case $\rho > 0$ Descartes' rule of signs tells us that the numerator of (3.1) has 2, 4, or 6 real zeros, exactly one of which is negative. Since $z=1$ and $z=-1$ are always zeros of this numerator, $p_{2n}(\rho,z)$ has no real zeros at all, or two or four real zeros, irrespective of how large $n$ is. For $\rho > 0$, the connections between the few zeros that are real to the eigenvalues of $K_n(\rho)$ will be explored in Section 6. But a zero that *lies on the unit circle* has a straightforward interpretation provided only that $\rho \neq \pm 1$:

**Proposition 3.4** Let $\rho \in \mathbf{C} \setminus \{-1,1\}$, let $z_k$ be a zero of $p_{2n}(\rho,z)$, and let $\lambda_k$ be the corresponding eigenvalue given in Theorem 3.3. Then $|z_k| = 1$ iff $\lambda_k \in \text{range}\{\sigma(\rho,\theta)\}$ (see Remark 2.2).

*Proof:* For nonzero $\rho$ the desired equivalence follows by writing the denominator of the first expression in (3.12) as $1 - \rho(z_k + z_k^{-1}) + \rho^2$, setting $z_k = \exp(i\mu_k)$, and comparing with (2.8). By (1.1), (3.1), and (2.8), the equivalence also holds for the trivial case $\rho = 0$. □

Theorem 2.3(v) is a necessary and sufficient condition for the positive definiteness of $K_n(\rho)$. With Cohn's theorem, we obtain more such conditions:

**Theorem 3.5** Let $\rho \in \mathbf{C} \setminus \{-1,1\}$. The following four statements are equivalent:

(i) All eigenvalues $\lambda_0, \lambda_1, \ldots, \lambda_{n-1}$ of $K_n(\rho)$ lie on $\text{range}\{\sigma(\rho,\theta)\}$.

(ii) All zeros $z_0, z_1, \ldots, z_{2(n-1)}$ of $p_{2n}(\rho,z)$ lie on the unit circle.

(iii) $-1 < \rho < 1$.

(iv) $K_n(\rho)$ is positive definite.

*Proof:* (ii) $\Rightarrow$ (iii): If (ii) holds, then Cohn's theorem implies that (a) $p_{2n}(\rho,z)$ is self-inversive, and (b) all zeros of $dp_{2n}(\rho,z)/dz$ lie on the unit circle. By (a) and Lemma 3.1 we have $\rho \in \mathbf{R}$. Differentiating (3.1) we find the explicit expression

$$\frac{dp_{2n}(\rho,z)}{dz} = \frac{2(z-\rho)\left[1-\rho z + z^{2n}(\rho z - 1) + nz^{2n-1}(z^2-1)(z-\rho)\right]}{(z^2-1)^2},$$

which shows that $z = \rho \in \mathbf{R} \setminus \{-1,1\}$ is a zero of $dp_{2n}(\rho,z)/dz$. It follows from (b) that $-1 < \rho < 1$.

(iii) $\Rightarrow$ (ii): (The proof that follows resembles the proof of Proposition 3.1 of [28], but the latter proof concerns $s_{n+1}(\rho,z)$ rather than $p_{2n}(\rho,z)$.) We show that $-1<\rho<1$ and $p_{2n}(\rho,z)=0$ entail $|z|=1$. Eqn. (3.1) and $p_{2n}(\rho,z)=0$ imply

$$z^{2n}(z-\rho)^2 - (\rho z - 1)^2 = 0.$$

This equation cannot be satisfied if $z=\rho$ unless $\rho=\pm 1$, which contradicts our assumption $-1<\rho<1$. Thus $z\neq\rho$, so that

$$|z|^{2n} = \frac{|\rho z - 1|^2}{|z-\rho|^2}.$$

By transforming the right-hand side, it is easy to show that

$$|z|^{2n} - 1 = -\frac{(1-\rho^2)(|z|^2 - 1)}{|z-\rho|^2}. \tag{3.15}$$

The two sides of (3.15) can have the same sign only if $|z|=1$, which is what we wanted to show.

(iii) $\Leftrightarrow$ (iv): This equivalence has already been shown as Theorem 2.3(v).

(i) $\Leftrightarrow$ (ii): This equivalence is an immediate consequence of Proposition 3.4. □

**Remark 3.6** When $-1<\rho<1$, $K_n(\rho)$ is Hermitian, the corresponding Laurent matrix has a well defined and bounded symbol $\sigma(\rho,\theta)$, and range$\{\sigma(\rho,\theta)\}$ is a line segment lying on the real axis. A well-known theorem (see, e.g. [2], [20], or [21, p. 46]) states that the eigenvalues of all such matrices lie on range$\{\sigma(\rho,\theta)\}$. This provides another way of showing (iii) $\Rightarrow$ (i). The aforementioned theorem is a finite-dimensional analogue to the Brown-Halmos theorem for infinite-Toeplitz matrices [21, p. 21].

When $\rho\in\mathbf{R}$ with $|\rho|>1$, range$\{\sigma(\rho,\theta)\}$ is a segment of the real axis. By Theorem 3.5, at least one eigenvalue of $K_n(\rho)$ does not belong to this segment. As a matter of fact, in Section 6 we show that only one or two eigenvalues lie outside the segment. On the other hand, when $\rho \in \mathbf{C}$ with $|\rho| < 1$, the eigenvalues of $K_n(\rho)$ typically *do not* lie on $\text{range}\{\sigma(\rho,\theta)\}$, see the figures of [16].

From (3.2) we see that each zero $z_k$ of $p_{2n}(\rho,z)$ is a zero of the product $s_{n+1}(\rho,z)c_{n+1}(\rho,z)$, which has the two additional zeros $z_k = \pm 1$. The theorem that follows describes the allocation of the $z_k$'s between $s_{n+1}(\rho,z)$ and $c_{n+1}(\rho,z)$. (As shown in Section 4, this has implications for the eigenvectors.) Hereinafter, we denote

$$\xi_n = \frac{n+1}{n-1}. \tag{3.16}$$

We will also use the usual floor and ceiling notations, with

$$\lfloor n/2 \rfloor = \begin{cases} n/2, & n = \text{even} \\ (n-1)/2, & n = \text{odd} \end{cases} \qquad \lceil n/2 \rceil = \begin{cases} n/2, & n = \text{even} \\ (n+1)/2, & n = \text{odd} \end{cases}. \tag{3.17}$$

**Theorem 3.7** Let $\rho \in \mathbf{C} \setminus \{-\xi_n, -1, 1, \xi_n\}$ and let $z_0, z_1, \ldots, z_{2(n-1)}$ be the $2n$ zeros of $p_{2n}(\rho,z)$. Then $z_k \neq \pm 1$ for all $k$. Furthermore, $z_k$ is either a "type-1 zero" or a "type-2 zero." The two types, which are mutually exclusive, are defined as follows:

<u>Type 1:</u> $z_k$ satisfies $s_{n+1}(\rho, z_k) = s_{n+1}(\rho, z_k^{-1}) = 0$; in this case the $\lambda_k$ given by

$$\lambda_k = \frac{z_k^{1-n}(1-\rho^2)}{(z_k - \rho)^2} = \frac{z_k^{n+1}(1-\rho^2)}{(1-\rho z_k)^2} \tag{3.18}$$

is a "<u>type-1 eigenvalue</u>" of $K_n(\rho)$.

<u>Type 2:</u> $z_k$ satisfies $c_{n+1}(\rho, z_k) = c_{n+1}(\rho, z_k^{-1}) = 0$, in which case the $\lambda_k$ given by

$$\lambda_k = -\frac{z_k^{1-n}(1-\rho^2)}{(z_k - \rho)^2} = -\frac{z_k^{n+1}(1-\rho^2)}{(1-\rho z_k)^2} \tag{3.19}$$

is a "<u>type-2 eigenvalue</u>" of $K_n(\rho)$.

There are $2\lfloor n/2 \rfloor$ type-1 $z_k$'s and $2\lceil n/2 \rceil$ type-2 $z_k$'s, corresponding to $\lfloor n/2 \rfloor$ type-1 $\lambda_k$'s and $\lceil n/2 \rceil$ type-2 $\lambda_k$'s.

*Proof:* We have $p_{2n}(\rho, \pm 1) = 0$ only in the excluded cases $\rho = \pm 1$ and $\rho = \pm \xi_n$, as

$$p_{2n}(\rho, \pm 1) = (n-1)(\rho \mp 1)(\rho \mp \xi_n).$$

by (2.5). Thus $z_k \neq \pm 1$. Furthermore, (3.3) and (3.4) imply that $s_{n+1}(\rho, z)$ and $c_{n+1}(\rho, z)$ have no common zeros unless $\rho = \pm 1$. Hence types 1 and 2 are mutually exclusive. We prove what remains separately for even and odd $n$: If $n = 2, 4, \ldots$ then (3.3) and (3.4) yield $s_{n+1}(\rho, 1) = c_{n+1}(\rho, -1) = 0$. We thus rearrange (3.2) according to

$$p_{2n}(\rho, z) = \frac{s_{n+1}(\rho, z)}{z-1} \frac{c_{n+1}(\rho, z)}{z+1}. \tag{3.20}$$

where the two factors are $n$'th degree polynomials. Since $z_k \neq \pm 1$, neither factor can have a zero equal to $\pm 1$. Therefore the $n$ zeros of the first factor $\frac{s_{n+1}(\rho, z)}{z-1}$ coincide with $n$ zeros of $s_{n+1}(\rho, z)$ (specifically, with those zeros of $s_{n+1}(\rho, z)$ that differ from 1). Similarly, the $n$ zeros of $\frac{c_{n+1}(\rho, z)}{z+1}$ are zeros of $c_{n+1}(\rho, z)$. We have thus shown the desired allocation for $n$ =even. Eqn. (3.7) then implies that the zeros thus allocated occur in inverse pairs.

For $n = 3, 5, \ldots$ we have $s_{n+1}(\rho, \pm 1) = 0$. Thus in place of (3.20) we write

$$p_{2n}(\rho, z) = \frac{s_{n+1}(\rho, z)}{z^2 - 1} c_{n+1}(\rho, z), \tag{3.21}$$

which shows that the zeros of $p_{2n}(\rho, z)$ consist of: (a) the $n - 1 = 2\lfloor n/2 \rfloor$ zeros of $s_{n+1}(\rho, z)$ that differ from $\pm 1$; and (b) the $n + 1 = 2\lceil n/2 \rceil$ zeros of $c_{n+1}(\rho, z)$.

If $s_{n+1}(\rho, z_k) = 0$ then (3.3) gives $1 - \rho z_k = z_k^n(z_k - \rho)$. Combining this with (3.12) gives (3.18). Eqn. (3.19) is similarly shown from (3.4) and (3.12). □

## 4. Trigonometric functions; eigensystem; double eigenvalues

Theorem 4.1 below discusses *complex* zeros $\mu_k$ of $s_{n+1}(\rho, e^{i\mu})$ and $c_{n+1}(\rho, e^{i\mu})$ (the significance of *real* zeros $\mu_k$ was noted in Proposition 3.4). By (3.3) and (3.4),

$$s_{n+1}(\rho, e^{i\mu}) = 2i \exp\left(i\frac{n+1}{2}\mu\right)\left[\sin\frac{\mu(n+1)}{2} - \rho \sin\frac{\mu(n-1)}{2}\right], \quad (4.1)$$

$$c_{n+1}(\rho, e^{i\mu}) = 2 \exp\left(i\frac{n+1}{2}\mu\right)\left[\cos\frac{\mu(n+1)}{2} - \rho \cos\frac{\mu(n-1)}{2}\right], \quad (4.2)$$

so we are led to examine the roots of the trigonometric functions in square brackets. For the special case $0 < \rho < 1$, the relations of these functions to the eigensystem have been discussed in many works [1]-[4], [29] (see also the slightly different functions of signal-processing works such as [5]-[7]). Theorem 4.1—which uses the terms introduced in Theorem 3.7—thus demonstrates that the well-known trigonometric functions pertain, more generally, to the eigensystem for all $\rho \in \mathbf{C}$.

**Theorem 4.1** Let $\rho \in \mathbf{C} \setminus \{-\xi_n, -1, 1, \xi_n\}$, $z_k \in \mathbf{C} \setminus \{0\}$, and let $\mu_k \in \mathbf{C}$ be defined by

$$\mu_k = -i \ln(z_k) \Leftrightarrow z_k = \exp(i\mu_k). \quad (4.3)$$

Then

(i) $z_k$ is a <u>type-1 zero</u> of $p_{2n}(\rho, z)$ iff $z_k \neq \pm 1$ ($\mu_k \neq 0, \pm\pi, \pm 2\pi, \ldots$) and

$$\sin\frac{\mu_k(n+1)}{2} - \rho \sin\frac{\mu_k(n-1)}{2} = 0. \quad (4.4)$$

In this case, the corresponding <u>type-1 eigenvalue</u> $\lambda_k$ can be found from

$$\lambda_k = \frac{1-\rho^2}{1-2\rho\cos\mu_k+\rho^2} = -\frac{\sin(n\mu_k)}{\sin\mu_k}, \quad (4.5)$$

and the vector $\mathbf{y}_k$ whose elements $y_{kj}$ are given by

$$y_{kj} = \sin\left[\mu_k\left(j - \frac{n-1}{2}\right)\right], \quad j = 0, 1, \ldots, n-1 \quad (4.6)$$

is a $\lambda_k$-eigenvector.

(ii) $z_k$ is a <u>type-2 zero</u> of $p_{2n}(\rho, z)$ iff $z_k \neq \pm 1$ ($\mu_k \neq 0, \pm\pi, \pm 2\pi, \ldots$) and

$$\cos\frac{\mu_k(n+1)}{2} - \rho\cos\frac{\mu_k(n-1)}{2} = 0. \tag{4.7}$$

In this case, the corresponding <u>type-2 eigenvalue</u> $\lambda_k$ can be found from

$$\lambda_k = \frac{1-\rho^2}{1-2\rho\cos\mu_k + \rho^2} = \frac{\sin(n\mu_k)}{\sin\mu_k}. \tag{4.8}$$

and the vector $\mathbf{y}_k$ whose elements $y_{kj}$ are given by

$$y_{kj} = \cos\left[\mu_k\left(j - \frac{n-1}{2}\right)\right], \quad j = 0,1,\ldots,n-1 \tag{4.9}$$

is a $\lambda_k$-eigenvector.

*Proof:* By (4.1) and Theorem 3.7, the two conditions $z_k \neq \pm 1$ and (4.4) are equivalent to $p_{2n}(\rho, z_k) = 0$ and $s_{n+1}(\rho, z_k) = 0$. By Theorem 3.3, the $\lambda_k$ in (3.12) is well defined and is an eigenvalue; with $z_k = \exp(i\mu_k)$, this $\lambda_k$ gives the first expression in (4.5). Eqn. (4.4) gives

$$\rho = \sin\frac{\mu_k(n+1)}{2} \Big/ \sin\frac{\mu_k(n-1)}{2}. \tag{4.10}$$

Substitution of (4.10) into the first expression (4.5) and use of trigonometric identities leads to the second expression in (4.5).

The assumption $\rho \neq \pm 1$ and Theorem 2.1(iv) show that $K_n(\rho)$ is nonsingular. Consequently, to show $\left[K_n(\rho) - \lambda_k I_n\right]\mathbf{y}_k = 0$ it suffices to prove

$$\left(\left[K_n(\rho)\right]^{-1} - \lambda_k^{-1} I_n\right)\mathbf{y}_k = 0. \tag{4.11}$$

The first formula for $\lambda_k$ in (4.5), $\rho \neq \pm 1$, and (2.1) show that (4.11) is tantamount to

$$\begin{bmatrix} 2\rho\cos\mu_k - \rho^2 & -\rho & 0 & \cdots & 0 & 0 & 0 \\ -\rho & 2\rho\cos\mu_k & -\rho & \cdots & 0 & 0 & 0 \\ 0 & -\rho & 2\rho\cos\mu_k & \cdots & 0 & 0 & 0 \\ \vdots & \vdots & \vdots & \ddots & \vdots & \vdots & \vdots \\ 0 & 0 & 0 & \cdots & 2\rho\cos\mu_k & -\rho & 0 \\ 0 & 0 & 0 & \cdots & -\rho & 2\rho\cos\mu_k & -\rho \\ 0 & 0 & 0 & \cdots & 0 & -\rho & 2\rho\cos\mu_k - \rho^2 \end{bmatrix} \mathbf{y}_k = 0. \tag{4.12}$$

It remains to show the $n$ equations in (4.12). Upon substituting $\mathbf{y}_k$ and $\rho$ by their values given by (4.6) and (4.10) it is readily understood that the desired $n$ equations are trigonometric identities (which hold for arbitrary $\mu_k$ and $n$). This completes our proof pertaining to type 1. The proof for type 2 is completely analogous. □

**Remark 4.2** The vector $\mathbf{y} = [y_0, y_1, \ldots, y_{n-1}]^T \in \mathbf{C}^n$ is called symmetric or even if $y_{n-1-j} = y_j$ and skew-symmetric or odd if $y_{n-1-j} = -y_j$. Thus type-1 (type-2) eigenvalues correspond to odd (even) $\mathbf{y}_k$'s. Accordingly, for $\rho \in \mathbf{C} \setminus \{0, \pm 1, \pm \xi_n\}$ we will distinguish type-1 and type-2 zeros and eigenvalues by using an odd and even index $k$, respectively (recall from Theorem 3.7 that the two types are mutually exclusive). In Section 6, this notation will prove to be useful for the case $\rho \in \mathbf{R}$.

**Remark 4.3** By Theorems 3.7 and 4.1, there are $\lfloor n/2 \rfloor$ $\lambda_k$ whose $\mathbf{y}_k$ is skew-symmetric, and $\lceil n/2 \rceil$ $\lambda_k$ whose $\mathbf{y}_k$ is symmetric. For the special case $\rho \in \mathbf{R}$ this was expected because $K_n(\rho)$ is real symmetric and centrosymmetric [30].

**Remark 4.4** For the special case $0 < \rho < 1$, our type-1/type-2 eigenvalues are the "odd/even" eigenvalues of Trench [3]. The terms odd/even are used, more generally, for eigenvalues of all *real* symmetric Toeplitz matrices [31].

It is easy to verify that $K_3(i2\sqrt{2})$ and $K_4(-1+2i)$ have the double eigenvalues $\lambda = -3$ and $\lambda = -4$, respectively. These facts exemplify the following

**Theorem 4.5** Let $\rho \in \mathbf{C} \setminus \{-1, 0, 1\}$. If $\lambda$ is a repeated eigenvalue of $K_n(\rho)$, then

$$\lambda = -n \qquad (4.13)$$

and $\lambda$ is a double eigenvalue. Moreover, $\lambda = -n$ is a double eigenvalue of <u>type 1</u> iff

$$\rho = \xi_n \frac{T_{(n+1)/2}(t_0)}{T_{(n-1)/2}(t_0)}, \qquad (4.14)$$

where $\xi_n$ is given in (3.16) and where $t_0 \in \mathbf{C}$ is any zero of the polynomial

$$q_1(n,t) = \begin{cases} \dfrac{U_{n-1}(t)-n}{t-1}, & n = 4,6,8\ldots \\ \dfrac{U_{n-1}(t)-n}{t^2-1}, & n = 5,7,9\ldots; \end{cases} \tag{4.15}$$

while $\lambda = -n$ is a double eigenvalue of <u>type 2</u> iff

$$\rho = \dfrac{T_{(n+1)/2}(t_0)}{T_{(n-1)/2}(t_0)}, \tag{4.16}$$

where $t_0 \in \mathbf{C}$ is any zero of the polynomial

$$q_2(n,t) = \begin{cases} \dfrac{U_{n-1}(t)+n}{t+1}, & n = 4,6,8\ldots \\ U_{n-1}(t)+n, & n = 3,5,7,\ldots. \end{cases} \tag{4.17}$$

In (4.14) and (4.16) we denote

$$T_{(n\pm 1)/2}(t) = \cos\left(\dfrac{n\pm 1}{2}\operatorname{Arccos} t\right), \quad t \in \mathbf{C}. \tag{4.18}$$

The multivaluedness of $\operatorname{Arccos} t$ does not affect (4.14) and (4.16). When $n = $ odd, $T_{(n\pm 1)/2}(t)$ is the Chebyshev polynomial of the first kind.

*Proof:* Let $\rho \neq 0, \pm 1$ and let $\lambda$ be a repeated eigenvalue. By Proposition 6.1 below, $K_n(\pm \xi_n)$ does not have repeated eigenvalues, so $\rho \neq \pm \xi_n$. By Theorem 3.7 the eigenvalue $\lambda$ is either of type 1 or of type 2. By Theorem 4.1, the repeated eigenvalue $\lambda$ is of type 1 iff

$$\lambda = -\dfrac{\sin(n\mu)}{\sin \mu}, \tag{4.19}$$

where $\mu \neq 0, \pm\pi,\ldots$ satisfies both (4.4) and its derivative so that

$$\sin\dfrac{\mu(n+1)}{2} = \rho \sin\dfrac{\mu(n-1)}{2}, \tag{4.20}$$

$$\xi_n \cos\dfrac{\mu(n+1)}{2} = \rho \cos\dfrac{\mu(n-1)}{2}. \tag{4.21}$$

Without loss of generality we assume $\Re\mu \in (0, \pi)$. If the multiplicity of $\lambda$ were larger than 2, $\mu$ would concurrently satisfy (4.20) and the second derivative of (4.20), something impossible. Thus $\lambda$ is a double eigenvalue. Set $\mu = \text{Arccos}\, t_0$ and solve (4.21) for $\rho$ to obtain (4.14) with (4.18). Letting $\chi = \tan\left[\frac{(n-1)\mu}{2}\right] / \tan\left[\frac{(n+1)\mu}{2}\right]$ allows us to rewrite (4.19) as $\lambda = (\chi+1)/(\chi-1)$ and to obtain $\chi = 1/\xi_n$ from (4.20) and (4.21). Thus $\lambda = (1+\xi_n)/(1-\xi_n)$ and (3.16) gives (4.13). Eqns. (4.13), (4.19), (2.6), and $t_0 = \cos\mu$ give $U_{n-1}(t_0) - n = 0$. This polynomial equation is satisfied by $t_0 = 1$ (or $\mu = 0$) when $n = $ even and by $t_0 = \pm 1$ (or $\mu = 0, \pi$) when $n = $ odd. Because we have assumed $\Re\mu \in (0, \pi)$, we exclude the aforesaid values of $t_0$. Therefore $t_0$ is a zero of $q_1(n,t)$, the polynomial given in (4.15). We have thus shown all assertions pertaining to type-1 eigenvalues. The proof for the type-2 case is entirely similar. □

## 5. Large eigenvalues for $|\rho| > 1$ and $n \to \infty$

Let $\rho \in \mathbb{C}$ with $|\rho| > 1$, and let $z \in \mathbb{C}$ be independent of $n$ with $|z| < 1$. As $n \to \infty$, the first two terms in (3.3) and (3.4) are exponentially smaller than the last two. Therefore

$$\begin{Bmatrix} s_{n+1}(\rho,z) \\ c_{n+1}(\rho,z) \end{Bmatrix} \sim \pm(\rho z - 1), \quad \text{as} \quad n \to \infty. \tag{5.1}$$

The right-hand side of this asymptotic approximation vanishes when $z$ assumes the value $1/\rho$ and, consistent with our initial assumptions, this value is independent of $n$ and lies within the unit circle (i.e., $|z| = |1/\rho| < 1$). In other words, the value $1/\rho$ is an asymptotic (large-$n$) zero of both $s_{n+1}(\rho,z)$ and $c_{n+1}(\rho,z)$. Substituting $z_k = z_1 \sim 1/\rho$ into the *first* expression (3.18) (the odd index $k = 1$ is consistent with Remark 4.2) gives an asymptotic expression for a type-1 eigenvalue of $K_n(\rho)$:

$$\lambda_1 \sim -\frac{\rho^{n+1}}{\rho^2 - 1}, \quad \text{as} \quad n \to \infty \quad (\rho \in \mathbf{C} \text{ with } |\rho| > 1). \tag{5.2}$$

This eigenvalue is negative and exponentially large. Substituting $z_k = z_1 \sim 1/\rho$ into (4.3) and (4.6) gives a large-$n$ approximation to the $\lambda_1$-eigenvector $\mathbf{y}_1$:

$$y_{1j} \sim \rho^{\frac{n-1}{2} - j} - \rho^{j - \frac{n-1}{2}}, \quad j = 0, 1, \ldots, n-1, \tag{5.3}$$

where the factor $1/(2i)$ was omitted. Similarly, substitution of $z_k = z_0 \sim 1/\rho$ into the *first* expression (3.19) leads to an asymptotic expression for a type-2 eigenvalue,

$$\lambda_0 \sim \frac{\rho^{n+1}}{\rho^2 - 1}, \quad \text{as} \quad n \to \infty \quad (\rho \in \mathbf{C} \text{ with } |\rho| > 1). \tag{5.4}$$

This positive and exponentially large eigenvalue corresponds to an eigenvector $\mathbf{y}_0$ whose large-$n$ approximation is found using (4.3) and (4.9):

$$y_{0j} \sim \rho^{\frac{n-1}{2} - j} + \rho^{j - \frac{n-1}{2}}, \quad j = 0, 1, \ldots, n-1. \tag{5.5}$$

where we omitted the factor $1/2$. Numerically, (5.2) and (5.4) give estimates that are very close to the two largest eigenvalues determined from standard routines: When $|\rho| = 3$ and $n = 10$, for example, relative errors in both eigenvalues' magnitudes remain less than 0.06% irrespective of the phase of $\rho$, and the value 0.06% decreases rapidly as $n$ increases. We also found that (5.3) and (5.5) can give accurate numerical estimates of the two corresponding eigenvectors.

**Remark 5.1** The subscripts 1 and 0 in (5.2)-(5.5) have been chosen judiciously, to be consistent with Section 6.

**Remark 5.2** Application of perturbation methods to the equations $s_{n+1}(\rho, z) = 0$ and $c_{n+1}(\rho, z) = 0$ leads to correction terms for $z_1$, $z_0$ which, in turn, lead to corrections for the eigenpairs (5.2)-(5.5). More details are beyond the scope of the present paper.

**Remark 5.3** If we initially assume $|z| > 1$ instead of $|z| < 1$ then the right-hand side of (5.1) is replaced by $z^{n+1} - \rho z^n$. This give the large-$n$ zero $z = \rho$ which, consistent with Theorem 3.3, is the reciprocal of the one we found previously. It leads to (5.2) and (5.4) by means of the *second* expressions in (3.18) and (3.19).

**Remark 5.4** In his work [18], Trench discusses the spectrum of the so-called generalized Kac-Murdoch-Szegö matrices. These form a class of *Hermitian* Toeplitz matrices containing (among other parameters) a real parameter $\phi$. When $\phi = 0$ this class reduces to a narrower class of matrices, all of which are real symmetric. This narrower class contains our $K_n(\rho)$ provided that $\rho \in \mathbf{R}$. Our (5.2) and (5.4) are consistent with a formula in [18] and, specifically, with the "sharpened" conclusion that can be found on p. 264 of [18].[1] In other words, our (5.2) and (5.4) are known for the special case $\rho \in \mathbf{R}$ with $|\rho| > 1$. On the other hand, Trench's proof does not pertain to the non-Hermitian case $\rho \in \mathbf{C} \setminus \mathbf{R}$, $|\rho| > 1$.

## 6. Real-symmetric case

The remainder of this paper deals with the real-symmetric (and Hermitian) case, which is the special case $\rho \in \mathbf{R}$. Oscillatory matrices have distinct real eigenvalues [26] so, by Theorem 2.3(ix), this is true of $K_n(\rho)$ if $0 < \rho < 1$ [3]. More generally,

**Proposition 6.1** For $\rho \in \mathbf{R} \setminus \{-1, 0, 1\}$, $K_n(\rho)$ has $n$ distinct real eigenvalues.

*Proof:* Eqn. (2.1) implies that $A = (1 - \rho^2)[K_n(\rho)]^{-1}$ is a real tridiagonal matrix whose elements $A_{jk}$ satisfy $A_{j,j+1} A_{j+1,j} = (-\rho)^2 > 0$. By [24, 3.1.P22, p. 174] or [20], such matrices (often called Jacobi matrices ) have distinct real eigenvalues. By $\rho \neq \pm 1$ and Theorem 2.1(iv), this property also holds true for $K_n(\rho)$. □

For $\rho \in \mathbf{R}$, we now give names for the eigenvalues encountered in Proposition 3.4.

**Definition 6.2** For $\rho \in \mathbf{R} \setminus \{-1, 1\}$, an eigenvalue of $K_n(\rho)$ is called *ordinary* if it lies within $\text{range}\{\sigma(\rho, \theta)\}$ (see Remark 2.2). In the trivial cases $\rho = 0, \pm 1$, all eigenvalues are called ordinary. An eigenvalue is called *extraordinary* if it is not ordinary.

---

[1] Note, however, that there is an obvious misprint in Trench's formula. When corrected, Trench's formula becomes consistent with ours.

We assume $\rho \geq 0$ because all results on the eigensystem of $K_n(\rho)$ carry over to $\rho < 0$ by means of Lemma 2.5. It is a consequence of (2.8) and Proposition 3.4 that

**Proposition 6.3** Let $\rho \geq 0$, let $\lambda_k \in \mathbf{R}$ ($k = 0, 1, \ldots, n-1$) be an eigenvalue of $K_n(\rho)$ corresponding to the zero $z_k$ of $p_{2n}(\rho, z)$, and let $\mu_k$ be given by (4.3). Then the following three statements are equivalent:

(i) $\lambda_k$ is an ordinary eigenvalue;

(ii) $\min\left\{\dfrac{1-\rho}{1+\rho}, \dfrac{1+\rho}{1-\rho}\right\} \leq \lambda_k \leq \max\left\{\dfrac{1-\rho}{1+\rho}, \dfrac{1+\rho}{1-\rho}\right\}$ or $\rho = 1$; (6.1)

(iii) $\mu_k \in \mathbf{R}$, $|z_k| = 1$.

**Remark 6.4** Our "ordinary/extraordinary eigenvalues" roughly correspond to the distributed/outlying spectrum introduced by Trench in [18] for the "generalized Kac-Murdoch-Szegö matrices" that we discussed in Remark 5.4 (see also [17], as well as the notes in [19, p. 275]). The basic difference is that the defining property of Trench's outlying spectrum holds for sufficiently large matrix dimension. By contrast, $n$ is unrestricted in Definition 6.2. (Indeed, with the exception of the large-$n$ formulas, all results in this paper hold for arbitrary $n$, down to $n = 2$.)

For $\rho \geq 0$, we will now characterize each $\mu_k$ of Theorem 4.1 (where we more generally assumed $\rho \in \mathbf{C}$) as the *unique* solution to a transcendental equation from within a *pre-specified* interval. The various intervals will lie on the real or the imaginary axis. These results will allow us to obtain simple conditions (on $\rho$ and $n$, see Theorem 6.6 and Corollary 6.9) for $\lambda_k$ to be extraordinary. Also, the various formulas can be numerically implemented using root-finding techniques, especially bracketing methods. This section's central result is Theorem 6.5; it involves the trigonometric functions

$$c_n^{(t)}(\mu) = \frac{\cos\dfrac{\mu(n+1)}{2}}{\cos\dfrac{\mu(n-1)}{2}}, \qquad s_n^{(t)}(\mu) = \frac{\sin\dfrac{\mu(n+1)}{2}}{\sin\dfrac{\mu(n-1)}{2}}, \qquad (6.2)$$

as well as the hyperbolic functions $c_n^{(h)}(x) = c_n^{(t)}(ix)$ and $s_n^{(h)}(x) = s_n^{(t)}(ix)$, viz.,

$$c_n^{(h)}(x) = \frac{\cosh\frac{x(n+1)}{2}}{\cosh\frac{x(n-1)}{2}}, \qquad s_n^{(h)}(x) = \frac{\sinh\frac{x(n+1)}{2}}{\sinh\frac{x(n-1)}{2}}. \tag{6.3}$$

Theorem 6.5 further involves the $n$-dependent numbers $\alpha_k, \beta_k$, and $\gamma_k$ defined by

$$\alpha_k = \frac{(k-1)\pi}{n-1}, \quad k = 1, 2, 3, \ldots, n-1, \tag{6.4}$$

$$\beta_k = \frac{k\pi}{n}, \quad k = 0, \ldots, n-1, \tag{6.5}$$

$$\gamma_k = \frac{(k+1)\pi}{n+1}, \quad k = 0, 1, \ldots, n-1. \tag{6.6}$$

Evidently,

$$0 = \alpha_1 < \beta_1 < \alpha_2 < \beta_2 < \alpha_3 < \ldots < \beta_{n-2} < \alpha_{n-1} < \beta_{n-1} < \pi, \tag{6.7}$$

$$0 = \beta_0 < \gamma_0 < \beta_1 < \gamma_1 < \beta_2 \ldots < \gamma_{n-2} < \beta_{n-1} < \gamma_{n-1} < \pi. \tag{6.8}$$

Furthermore, from (6.2), (6.3), and (3.16) we see that

$$c_n^{(t)}(\alpha_k + 0) = +\infty \text{ for } k = \text{even}; \quad s_n^{(t)}(\alpha_k + 0) = +\infty \text{ for } k = \text{odd with } k \geq 3, \tag{6.9}$$

$$c_n^{(t)}(\beta_k) = 1 \text{ for } k = \text{even}; \quad s_n^{(t)}(\beta_k) = 1 \text{ for } k = \text{odd}, \tag{6.10}$$

$$c_n^{(t)}(\gamma_k) = 0 \text{ for } k = \text{even}; \quad s_n^{(t)}(\gamma_k) = 0 \text{ for } k = \text{odd}, \tag{6.11}$$

$$c_n^{(t)}(0) = c_n^{(h)}(0) = 1; \quad s_n^{(t)}(0) = s_n^{(h)}(0) = \xi_n, \tag{6.12}$$

$$c_n^{(h)}(+\infty) = s_n^{(h)}(+\infty) = +\infty. \tag{6.13}$$

Our exposition of Theorems 6.5 and 6.6 emphasizes the continuity of each $\lambda_k = \lambda_k(\rho)$ as function of the nonnegative variable $\rho$. To this end, Theorem 6.5 includes some known results for the case $0 < \rho < 1$ that can be found in Trench's work [3]. The said results coincide, specifically, with the "slightly improved estimates" of p. 4 of [3] or [10, eqn. (16)], which are refinements of results of Kac, Murdock, and Szegö [1]; the latter results (of [1]) can also be found in [2], [4], [29], and other works. To the best of the author's knowledge, the parts of Theorem 6.5 pertaining to the case $\rho > 1$ are new.

**Theorem 6.5** Let $\rho \geq 0$, and let "tur" stand for "the unique root." The $n$ eigenvalues $\lambda_k(\rho)$ of $K_n(\rho)$ can be found from

$$\lambda_k(\rho) = \frac{1-\rho^2}{1-2\rho\cos[\mu_k(\rho)]+\rho^2}, \quad k=0,1,\ldots,n-1, \tag{6.14}$$

or, alternatively, from

$$\lambda_k(\rho) = (-1)^k \frac{\sin[n\mu_k(\rho)]}{\sin[\mu_k(\rho)]}, \quad k=0,1,\ldots,n-1, \tag{6.15}$$

where the first two parameters $\mu_0(\rho)$ and $\mu_1(\rho)$ are given by

$$\mu_0(\rho) = \begin{cases} ix_0(\rho), \text{ where } x_0(\rho) \text{ is tur of } c_n^{(h)}(x)=\rho \text{ in } [0,+\infty), & \rho \geq 1 \\ \text{tur of } c_n^{(t)}(\mu)=\rho \text{ in } [0,\gamma_0], & 1 \geq \rho \geq 0, \end{cases} \tag{6.16}$$

$$\mu_1(\rho) = \begin{cases} ix_1(\rho), \text{ where } x_1(\rho) \text{ is tur of } s_n^{(h)}(x)=\rho \text{ in } [0,+\infty), & \rho \geq \xi_n \\ \text{tur of } s_n^{(t)}(\mu)=\rho \text{ in } [0,\beta_1], & \xi_n \geq \rho \geq 1 \\ \text{tur of } s_n^{(t)}(\mu)=\rho \text{ in } [\beta_1,\gamma_1], & 1 \geq \rho \geq 0. \end{cases} \tag{6.17}$$

The remaining parameters $\mu_2(\rho),\ldots,\mu_{n-1}(\rho)$ have different expressions according to whether $k$ is even or odd. In the first case,

$$\mu_k(\rho) = \begin{cases} \text{tur of } c_n^{(t)}(\mu)=\rho \text{ in } (\alpha_k,\beta_k], & \rho \geq 1 \\ \text{tur of } c_n^{(t)}(\mu)=\rho \text{ in } [\beta_k,\gamma_k], & 1 \geq \rho \geq 0 \end{cases}, \quad k=2,4,\ldots,2\lceil n/2 \rceil-2, \tag{6.18}$$

and in the second,

$$\mu_k(\rho) = \begin{cases} \text{tur of } s_n^{(t)}(\mu)=\rho \text{ in } (\alpha_k,\beta_k], & \rho \geq 1 \\ \text{tur of } s_n^{(t)}(\mu)=\rho \text{ in } [\beta_k,\gamma_k], & 1 \geq \rho \geq 0 \end{cases}, \quad k=3,5,\ldots,2\lfloor n/2 \rfloor-1. \tag{6.19}$$

Finally, for $\rho > 0$ with $\rho \neq 1$ and $\rho \neq \xi_n$, the $\lambda_k$ given above are of type 1/type 2 when $k$ = odd/even, respectively (the two types are defined in Theorem 3.7).

Theorem 6.5 translates into a specific ordering of the $\lambda_k$, yields upper and lower bounds, and gives simple conditions for $\lambda_k$ to be ordinary or extraordinary:

**Theorem 6.6** The matrix $K_n(\rho)$ has exactly two extraordinary eigenvalues ($\lambda_0(\rho)$ and $\lambda_1(\rho)$) when $\rho > \xi_n$, exactly one ($\lambda_0(\rho)$) when $1 < \rho \leq \xi_n$, and none when $0 \leq \rho \leq 1$. Furthermore,

$$\lambda_{n-1} = \lambda_{n-2} = \ldots = \lambda_1 = \lambda_0 = 1, \quad \text{when } \rho = 0; \tag{6.20}$$

$$0 < \frac{1-\rho}{1+\rho} < \lambda_{n-1} < \lambda_{n-2} < \ldots < \lambda_1 < \lambda_0 < \frac{1+\rho}{1-\rho}, \quad \text{when } 0 < \rho < 1; \tag{6.21}$$

$$\lambda_1 = \lambda_2 = \ldots = \lambda_{n-1} = 0 < \lambda_0 = n, \quad \text{when } \rho = 1; \tag{6.22}$$

$$-\frac{\rho+1}{\rho-1} < \lambda_1 < \lambda_2 < \ldots < \lambda_{n-2} < \lambda_{n-1} < -\frac{\rho-1}{\rho+1} < 0 < n < \lambda_0, \quad \text{when } 1 < \rho < \xi_n; \tag{6.23}$$

$$-n = \lambda_1 < \lambda_2 < \lambda_3 < \ldots < \lambda_{n-2} < \lambda_{n-1} < -\frac{1}{n} < 0 < n < \lambda_0, \quad \text{when } \rho = \xi_n; \tag{6.24}$$

$$\lambda_1 < -n < -\frac{\rho+1}{\rho-1} < \lambda_2 < \lambda_3 < \ldots < \lambda_{n-2} < \lambda_{n-1} < -\frac{\rho-1}{\rho+1} < 0 < n < \lambda_0, \quad \text{when } \rho > \xi_n. \tag{6.25}$$

Before proving Theorems 6.5 and 6.6, let us add some comments that aid in understanding the theorems' notation and comprehending their results:

**Remarks 6.7**

(i) For $n = 5$, Fig. 1 is a plot of the eigenvalues $\lambda_0, \lambda_1, \ldots, \lambda_4$ as function of $\rho$ for $0 \leq \rho < 1.6$. In accordance with Theorem 6.6, the ordering of $\lambda_1, \lambda_2, \lambda_3, \lambda_4$ is inverted when $\rho = 1$. And, as the theorem predicts, $\lambda_1(0) = \ldots = \lambda_4(0) = 1$; $\lambda_0(1) = n = 5$; $\lambda_1(1) = \ldots = \lambda_4(1) = 0$; and $\lambda_1(\xi_n) = \lambda_1(\xi_5) = \lambda_1(1.5) = -n = -5$. In Fig. 1, $\lambda_0$ and $\lambda_1$ are extraordinary when $\rho > 1$ and $\rho > \xi_5 = 1.5$, respectively.

(ii) An *extraordinary* eigenvalue is given in Theorem 6.5 in terms of a real root ($x_0$ or $x_1$) of a *hyperbolic* function ($c_n^{(h)}$ or $s_n^{(h)}$). The corresponding parameter $\mu$ (i.e., $\mu_0 = ix_0$ or $\mu_1 = ix_1$) is purely imaginary and by (4.3), the corresponding zero of $p_{2n}(\rho, z)$ (i.e., $z_0$ or $z_1$) is real and positive.

Since zeros of $p_{2n}(\rho,z)$ come in inverse pairs, $p_{2n}(\rho,z)$ has none, two, or four real zeros, as observed in Section 3 via Descartes' rule of signs.

(iii) In Theorem 6.5, any *ordinary* eigenvalue $\lambda_k$ is found from a *real* parameter $\mu_k$ which is a root of a *trigonometric* function ($c_n^{(t)}$ or $s_n^{(t)}$). By Proposition 6.3, the corresponding root $z_k$ of $p_{2n}(\rho,z)$ lies on the unit circle. Any such $\mu_k(\rho)$ belongs to $[\alpha_k,\beta_k]$ when $\rho \geq 1$ or to $[\beta_k,\gamma_k]$ when $1 \geq \rho \geq 0$. Special cases (for the endpoint values of $\rho$) are:

$$\mu_k(0)=\gamma_k, \quad \mu_k(1)=\beta_k, \quad k=0,1,\ldots,n-1 \tag{6.26}$$

$$\mu_1(\xi_n)=0, \quad \mu_k(\infty)=\alpha_k, \quad k=2,3,\ldots,n-1 \tag{6.27}$$

(iv) Let $\rho \neq 0$, $\rho \neq 1$, and $\rho \neq \xi_n$. Then type-1 (type-2) $\lambda_k$ have odd (even) indices $k$ and involve the functions $s_n^{(t)}$ and $s_n^{(h)}$ ($c_n^{(t)}$ and $c_n^{(h)}$).

(v) Consider the equations involving $c_n^{(h)}$ and $s_n^{(h)}$ in (6.16) and (6.17), respectively. As $n \to \infty$, it is apparent from (6.3) that both equations possess the large-$n$ solution $x \sim \ln \rho$. Thus $z_0 \sim 1/\rho$ and $z_1 \sim 1/\rho$, corresponding to the large-$n$ eigenvalues $\lambda_0$ and $\lambda_1$ found in Section 5. Hence Section 5 used consistent subscripts ($k=0,1$) with those herein.

(vi) The case $\rho=1$ is special in several respects. Let $\rho$ vary from values smaller than 1 to values larger than 1, as in Fig. 1. When $\rho=1$ (so that $\mu_0=0$ and $z_0=1$), all eigenvalues except $\lambda_0(\rho)$ coalesce at the origin. The fact that these eigenvalues invert their ordering (compare (6.21) to (6.23)) is consistent with the $O(\rho-1)$ behavior suggested by the bottom Fig. 1 (and demonstrated in Section 7.2 below). At the same time, $\lambda_0(\rho)$ becomes extraordinary; $\sigma(\rho,\theta)$ no longer has meaning as the symbol of a Laurent matrix (see Remark 2.2); the "ordinary interval" in (6.1) changes abruptly; and the two conjugate, unimodular zeros of $p_{2n}(\rho,z)$ (namely, $z_0=e^{i\mu_0}$ and $z_0^{-1}=e^{-i\mu_0}$) coalesce to become the double zero $z_0=1$ (they then split into two real and inverse zeros). This means that $\rho=1$ is a *critical point* [32] of the algebraic function $p_{2n}(\rho,z)$.

(vii) When $\rho = \xi_n$, $p_{2n}(\rho, z)$ has the double root $z_1 = 1$ ($\mu_1 = 0$), signaling the appearance of the second extraordinary eigenvalue $\lambda_1(\rho)$.

*Proof of Theorem 6.5:* We initially assume $\rho \neq 0$ and $\rho \neq 1$. Existence of a solution to each of the nine equations written in (6.16)–(6.19) is a simple consequence of the continuity of the function in the pertinent interval and the function values at the interval endpoints, which are listed in (6.9)-(6.13). To give an example, when $\xi_n \geq \rho > 1$ the equation $s_n^{(t)}(\mu) = \rho$ (which is the second equation in (6.17)) possesses a solution $\mu = \mu_1$ in $[0, \beta_1)$ because $s_n^{(t)}(\mu)$ is continuous in $[0, \beta_1)$ by (6.2) and (6.5), with $s_n^{(t)}(0) = \xi_n$ by (6.12) and $s_n^{(t)}(\beta_1) = 1$ by (6.10).

Now further assume $\rho \neq \xi_n$. The nonzero $\mu_k$ defined by (6.16)-(6.19) satisfy (4.4) (when $k = $ odd) or (4.7) (when $k = $ even). Accordingly, (6.14) and (6.15) are restatements of (4.5) and (4.8) of Theorem 4.1, with the odd/even subscripts and the type-1/type-2 terminology consistent with Theorem 4.1 and Remarks 4.2 and 4.3.

In the case $0 < \rho < 1$, the $n$ $\mu_k \in (\beta_k, \gamma_k)$ are distinct and, by (6.8), lie in $(0, \pi)$; via (6.14), these $\mu_k$ lead to $n$ distinct eigenvalues $\lambda_k$. Arguing in a similar manner via (6.7), we see that our $\mu_k$ lead to $n$ distinct $\lambda_k$ for each of the two cases $1 < \rho < \xi_n$ and $\rho > \xi_n$ (note that $\mu_0$ is purely imaginary in both cases and that $\mu_1$ is purely imaginary in the second case). We have thus found the $n$ discrete and real $\lambda_k$ for all three cases. Finally, our $\mu_k$ are the *unique* solutions to the nine equations in (6.16)-(6.19) because any additional solution would mean more than $n$ eigenvalues.

Theorem 4.1 cannot be invoked when $\rho = \xi_n$ but, by (4.1), (4.2), and (3.2), the $n-1$ distinct and nonzero $\mu_0, \mu_2, \mu_3, \ldots, \mu_{n-1}$ still give $n-1$ distinct zeros $z_k = \exp(i\mu_k) \neq 1$ of $p_{2n}(\rho, z)$. Furthermore, the value $\mu_1 = 0$ of (6.17) corresponds to $z_1 = 1$ which, by (3.5), is yet another zero of $p_{2n}(\rho, z)$. For all $k$, we can use (3.12) (Theorem 3.3) to find $n$ distinct $\lambda_k$. Since (3.12) with $z_k = \exp(i\mu_k)$ reduces to (6.14), our theorem continues to hold.

In the trivial cases $\rho = 0$ and $\rho = 1$, substitution of (6.26) into (6.14) give $n$ (non-distinct) values $\lambda_k$, which are the true eigenvalues by (2.7). In particular, for the case $\mu_0(1) = \beta_0 = 0$, (6.14) becomes indeterminate; but (6.15), understood in a limiting sense, gives the correct value $\lambda_0(1) = n$. Thus our theorem still holds in the two trivial cases. □

*Proof of Theorem 6.6:* We only prove this theorem for $1 < \rho < \xi_n$ because the remaining cases are similar. Since $\mu_k \in (\alpha_k, \beta_k)$ for $k = 1, \ldots, n-1$, (6.7) gives
$$0 < \mu_1 < \mu_2 < \ldots < \mu_{n-1} < \pi.$$
As $\rho > 1$, (6.14) retains this ordering:
$$-(\rho+1)/(\rho-1) < \lambda_1 < \lambda_2 < \ldots < \lambda_{n-2} < \lambda_{n-1} < -(\rho-1)/(\rho+1) < 0.$$
By Proposition 6.3, all $\lambda_k$ appearing in the above inequality are ordinary. Eqn. (6.16) gives $\mu_0 = ix_0$ with $x_0 > 0$. Consequently (6.15) implies
$$\lambda_0 = \frac{\sinh(nx_0)}{\sinh x_0} > n,$$
which completes our proof of (6.23). Eqn. (6.23) and Proposition 6.3 now imply that $\lambda_0$ is the only extraordinary eigenvalue. □

When $\rho$ is close to 1, the rightmost (leftmost) bound in (6.21) (in (6.23)) can be large in absolute value. The following $n$-dependent bounds are often sharper:

**Lemma 6.8** In addition to the bounds given by Theorem 6.6, we have
$$\lambda_0 \leq n \quad \text{when} \quad 0 \leq \rho \leq 1, \tag{6.28}$$
$$\lambda_1 \geq -n \quad \text{when} \quad 1 \leq \rho \leq \xi_n. \tag{6.29}$$

*Proof:* Inequality (6.28) holds because $K_n(\rho)$ is positive semidefinite with trace equal to $n$. For $1 \leq \rho \leq \xi_n$, (6.17), (6.15), and (6.5) show that
$$\lambda_1 = -\frac{\sin(n\mu)}{\sin \mu}, \tag{6.30}$$
for some $\mu \in [0, \pi/n]$. Inequality (6.29) then follows from the easily-proven inequality
$$\frac{\sin(n\mu)}{\sin \mu} \leq n, \quad n = 2, 3, \ldots, \quad \mu \in \left[0, \frac{\pi}{n}\right]. \quad □$$

Theorem 6.6 and Lemma 6.8 now lead to a simple criterion—already noticeable in the example of Fig. 1—that distinguishes between ordinary and extraordinary eigenvalues:

**Corollary 6.9** Let $\rho \geq 0$. An eigenvalue $\lambda_k$ of $K_n(\rho)$ is extraordinary iff $|\lambda_k| > |n|$. In other words, the statement

(iv) $\quad |\lambda_k| \leq |n|$,

is equivalent to statements (i)-(iii) of Proposition 6.3.

*Proof:* By (6.20)-(6.25) and Lemma 6.8, $|\lambda_k| > |n|$ is true iff ($k = 0$ and $\rho > 1$) or ($k = 1$ and $\rho > \xi_n$). By Theorem 6.6, these are precisely the conditions for $\lambda_k$ to be extraordinary. Definition 6.2 now gives (iv) $\Leftrightarrow$ (i). □

For completeness, let us write down the eigenvectors for all known eigenvalues in the special cases $\rho = 1$ and $\rho = \xi_n$: For $\rho = 1$, $\mathbf{y}_0^T = [1,1,\ldots,1]$ is an eigenvector corresponding to $\lambda_0 = n$, while the eigenspace of the multiple eigenvalue $\lambda_1 = \ldots = \lambda_{n-1} = 0$ consists of all vectors whose elements sum to zero. For $\rho = \xi_n$, $\mathbf{y}_1^T = [-n+1, -n+3, \ldots, n-5, n-3, n-1]$ is an eigenvector corresponding to $\lambda_1 = -n$.

## 7. Real-symmetric case: Simple approximations for eigenvalues

We now build upon Theorem 6.5 to obtain some approximate formulas for the eigenvalues. They are simple enough to facilitate understanding and in some cases, precise enough for potential use as initial approximations for iterative eigenvalue solvers. The formulas that follow are exemplary—it is possible to obtain many more.

*7.1 Regula-falsi approximations to ordinary eigenvalues*

For the case $-1 < \rho < 1$, Trench [29] (see also [4]) proposed using the regula falsi method to the equations first developed by Kac, Murdock, and Szegö [1] as a means of *numerically* determining the eigenvalues. Here, for all $\rho \geq 0$ (for $\rho < 0$, use Lemma 2.5) we write down the *first approximation* that results when this iterative method is applied to the transcendental equations of Theorem 6.5. By (6.26) and (6.27), for the case of ordinary eigenvalues the said first approximations are

$$\mu_k(\rho) \cong \beta_k \rho + \gamma_k(1-\rho), \quad 0 \le \rho \le 1, \quad k = 0, 1, \ldots, n-1. \tag{7.1}$$

$$\mu_k(\rho) \cong \frac{\beta_k - \alpha_k}{\rho} + \alpha_k, \quad \rho \ge 1, \quad k = 2, 3, \ldots, n-1, \tag{7.2}$$

where $\alpha_k, \beta_k$, and $\gamma_k$ are given in (6.4)-(6.6). When combined with (6.14)—or with its alternative (6.15)—eqns. (7.1) and (7.2) constitute approximate formulas for $\lambda_k(\rho)$. As (7.1) and (7.2) result from intervals (namely, $[\beta_k, \gamma_k]$ and $[\alpha_k, \beta_k]$) whose lengths vanish as $n \to \infty$, our "approximate formulas" are, loosely speaking, large-$n$ formulas; but we make no attempt to render this statement more precise.

With (5.2) and (5.4), we have now obtained large-$n$ formulas for *all* $\lambda_k$, both ordinary and extraordinary.[2] Numerically, for $\rho = 3$ the estimates provided by (7.2) with (6.14) have a maximum error (over all $k$, as compared to the $\lambda_k$ obtained by standard routines) of 2.8% when $n = 10$, 0.71% when $n = 40$, and 0.18% when $n = 160$. For $\rho = 0.3$, the corresponding errors for (7.1) with (6.14) are 2.1% when $n = 10$, 0.55% when $n = 40$, and 0.14% when $n = 160$.

## 7.2 Approximations for $\rho$ close to 1

Small errors like those just reported occur when $\rho$ is far from 1. When $\rho \cong 1$ relative errors turned out larger, especially when $k$ is small. A formula appropriate for the case $\rho \to 1$ results by substituting $\mu_k = \beta_k = k\pi/n$ (which by (6.26) is the pertinent root when $\rho = 1$ and $k \ge 1$) into (6.14), expanding in powers of $\rho - 1$, and retaining the first term (this is legitimate because this first term is independent of higher-order terms in the expansion of $\mu_k$). We thus obtain the simple formula

$$\lambda_k \sim \frac{1-\rho}{1-\cos(k\pi/n)}, \quad k = 1, 2, \ldots, n-1, \quad \rho \to 1, \tag{7.3}$$

---

[2] The length of the interval in the second equation (6.17) is vanishingly small as $n \to \infty$. Consequently, that equation cannot produce a large-$n$ and fixed-$\rho$ formula similar to (7.1) and (7.2). For fixed $\rho \in \mathbf{R}$ we thus have two large-$n$ formulas for $\lambda_1$, namely (5.2) for the case $\rho > 1$, and (7.1) together with (6.14) for the case $0 < \rho < 1$.

which gives nearly identical results (at the scale of the figure) with the bottom Fig. 1. For the Perron root (case $k = 0$), set $\varepsilon = \rho - 1$ to obtain a small-$\varepsilon$ correction to $\lambda_0 = n$ which, by (6.22), is the eigenvalue when $\varepsilon = 0$. By (6.2) and (6.3),

$$c_n^{(h)}(x) = 1 + \frac{nx^2}{2} + O(x^4), \quad \text{as } x \to 0; \quad c_n^{(t)}(\mu) = 1 - \frac{n\mu^2}{2} + O(\mu^4), \quad \text{as } \mu \to 0.$$

Thus, both equations in (6.16) lead to the small-$\varepsilon$ root $\mu_0 = ix_0 = \sqrt{-2\varepsilon/n} + O(\varepsilon)$, which is purely imaginary (real) when $\varepsilon > 0$ ($\varepsilon < 0$). Substitution into (6.15) yields

$$\lambda_0(\rho) = n + \frac{n^2 - 1}{3}(\rho - 1) + O\left((\rho - 1)^{3/2}\right), \quad \rho \to 1. \quad (7.4)$$

When $n = 10$, (7.4) gives $\lambda_0(0.98) \cong 9.340$ and $\lambda_0(1.02) \cong 10.66$, values that present a 0.29% and 0.28% difference, respectively, to the numerically-determined Perron roots. Because our derivation assumed that $n$ is fixed, the numerical accuracy of (7.4) deteriorates if $n$ becomes too large.

## 8. Discussion

Our Definition 6.2 is only relevant to the *real symmetric* matrices $K_n(\rho)$, $\rho \in \mathbf{R}$. As mentioned in Section 3, in the general case $\rho \in \mathbf{C}$ eigenvalues typically do not lie on range$\{\sigma(\rho, \theta)\}$. Therefore the concept of extraordinary eigenvalues, as it stands, does not extend to the general case. Nonetheless, in the case $\rho \in \mathbf{C}$ with $|\rho| > 1$, and for sufficiently large matrix dimension $n$, numerical experiments indicated that:

(i) The simple criterion of Corollary 6.9 does not seem to carry over to the general case.

(ii) All eigenvalues lie *outside* of the convex hull of range$\{\sigma(\rho, \theta)\}$. (Contrast this to what happens in the case $|\rho| < 1$: In that case, the aforementioned (in Remark 3.6) finite-dimensional analogue of the Brown-Halmos theorem [20], [21, p. 46] predicts that all eigenvalues lie *within* the convex hull—note that all winding numbers are zero, so that the spectra of the infinite-Toeplitz and Laurent matrices coincide).

(iii) All eigenvalues except two are *close* to range$\{\sigma(\rho,\theta)\}$; as one would foresee, the two exceptions are closely approximated by (5.2) and (5.4).

The closeness mentioned in (iii) suggests that it might be possible to extend the notion of extraordinary eigenvalues to the case $\rho \in \mathbf{C}$ with $|\rho| > 1$. A satisfactory definition should lead to straightforward conditions (involving $\rho$, $n$, and the analytic continuation of the symbol) for the existence of exactly two extraordinary eigenvalues.

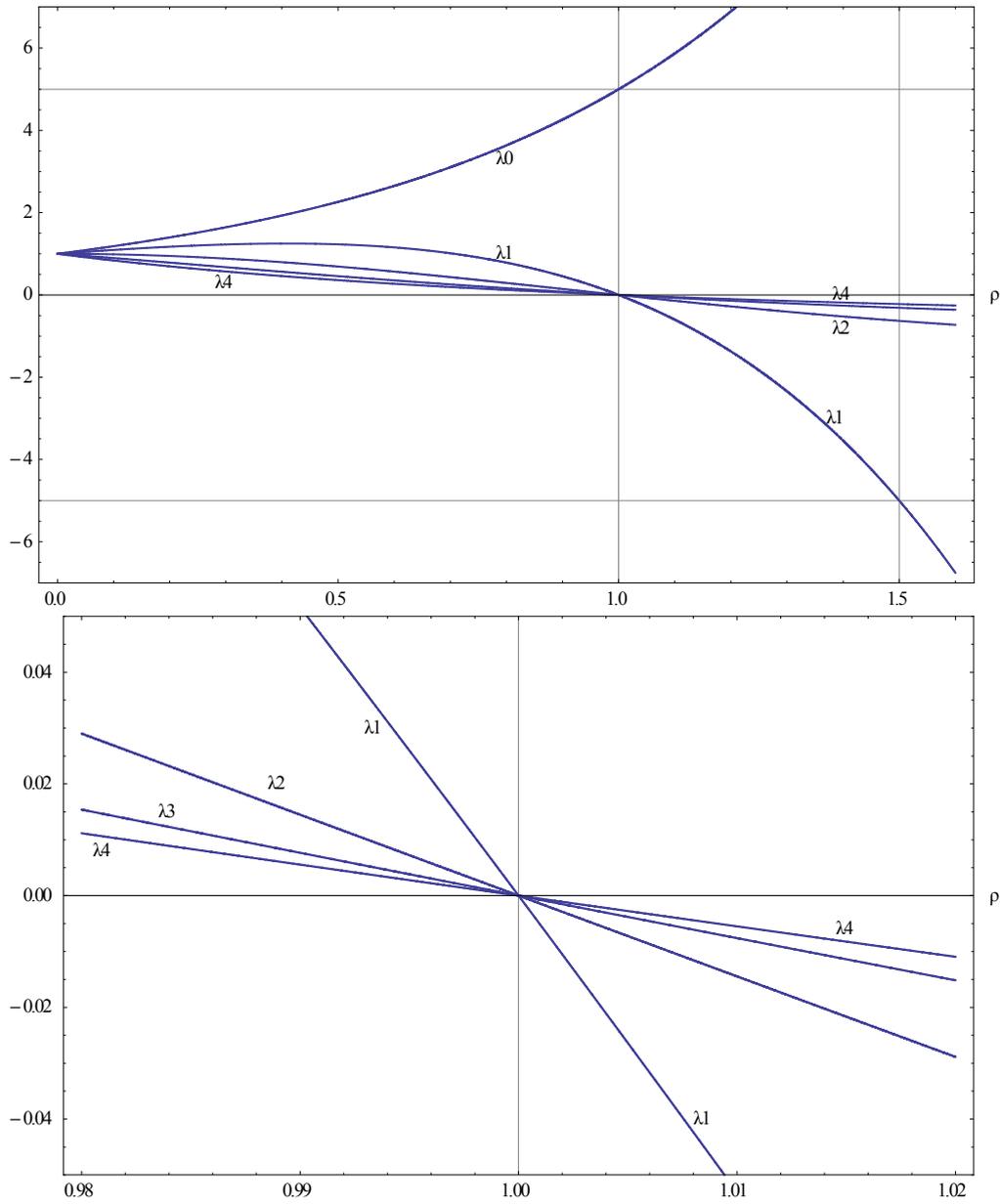

**Fig. 1** Eigenvalues $\lambda_0, \lambda_1, \lambda_2, \lambda_3, \lambda_4$ as function of $\rho$ for $n=5$. In the top figure, $0 \leq \rho < 1.6$. The bottom figure (where $\lambda_0$ is out of scale and not shown) is a closeup near $\rho = 1$.